\documentclass[11pt]{article}

\usepackage{latexsym,amssymb,amsmath,natbib,graphicx, amsthm, times}
\usepackage[margin=1in,left=1in]{geometry}
\usepackage{times}
\usepackage{natbib}
\renewcommand{\cite}{\citeyearpar}

\newcommand {\ME}{\mathbb{E}^{x}}

\renewcommand{\S}{\mathcal{S}}

\numberwithin{equation}{section}
\newtheorem{proposition}{Proposition}[section]
\newtheorem{coro}{Corollary}[section]
\newtheorem{remark}{Remark}[section]
\newtheorem{lemma}{Lemma}[section]
\newtheorem{example}{Example}[section]

\newcommand {\R}{\mathbb{R}}
\newcommand {\C}{\mathbf{C}}
\newcommand {\D}{e^{-\alpha\tau}}
\newcommand {\F}{\mathcal{F}}
\newcommand {\A}{\mathcal{A}}
\newcommand {\M}{\mathcal{M}}
\newcommand {\p}{\mathbb{P}}

\newcommand {\E}{\mathbb{E}}
\newcommand {\LL}{\mathcal{L}}

\title{A Direct Method for Solving Optimal Switching Problems of One-Dimensional Diffusions}
\author{
Masahiko Egami}
\date{}

\begin{document}

\maketitle
\begin{abstract}\noindent In this paper, we propose a direct solution method for optimal
switching problems of one-dimensional diffusions.  This method is
free from conjectures about the form of the value function and
switching strategies, or does not require the proof of optimality
through quasi-variational inequalities. The direct method uses a
general theory of optimal stopping problems for one-dimensional
diffusions and characterizes the value function as sets of the
smallest linear majorants in their respective transformed spaces.

\end{abstract}
\section{Introduction}  Stochastic \emph{optimal switching} problems (or \emph{starting and
stopping} problems) are important subjects both in mathematics and
economics. Since there are numerous articles about real options in
the economic and financial literature in recent years, the
importance and applicability of control problems including optimal
switching problems cannot be exaggerated.

A typical optimal switching problem is described as follows:   The
controller monitors the price of natural resources for optimizing
(in some sense) the operation of an extraction facility. She can
choose when to start extracting this resource and when to
temporarily stop doing so, based upon price fluctuations she
observes.  The problem is concerned with finding an optimal
switching policy and the corresponding value function.  A number of
papers on this topic are well worth mentioning : Brennan and Schwarz
\cite{BS1985} in conjunction with convenience yield in the energy
market, Dixit \cite{D1989} for production facility problems, Brekke
and {\O}ksendal \cite{BO1994} for resource extraction problems,
Yushkevich \cite{Y2001} for positive recurrent countable Markov
chain, and Duckworth and Zervos \cite{DZ2001} for reversible
investment problems.   Hamdad\`{e}ne and Jeanblanc \cite{HJ2004}
analyze a general adapted process for finite time horizon using
reflected stochastic backward differential equations. Carmona and
Ludkovski \cite{CL2005} apply to energy tolling agreement in a
finite time horizon using Monte-Carlo regressions.

A basic analytical tool for solving switching problems is
quasi-variational inequalities. This method is indirect in the sense
that one first conjectures the form of the value function and the
switching policy and next verifies the optimality of the candidate
function by proving that the candidate satisfies the variational
inequalities.  In finding the specific form of the candidate
function, appropriate boundary conditions including the smooth-fit
principle are employed.  This formation shall lead to a system of
non-linear equations that are often hard to solve and the existence
of the solution to the system is also difficult to prove.  Moreover,
this indirect solution method is specific to the underlying process
and reward/cost structure of the problem.  Hence a slight change in
the original problem often causes a complete overhaul in the highly
technical solution procedures.

Our solution method is direct in the sense that we first show a new
mathematical characterization of the value functions and, based on
the characterization, we shall \emph{directly} find the value
function and optimal switching policy.  Therefore, it is free from
any guesswork and applicable to a larger set of problems (where the
underlying process is one-dimensional diffusions) than the
conventional methods. Our approach here is similar to Dayanik and
Karatzas \cite{DK2003} and Dayanik and Egami \cite{DE2005} that
propose direct methods of solving optimal stopping problems and
stochastic impulse control problems, respectively.

The paper is organized in the following way.  In the next section,
after we introduce our setup of one dimensional optimal switching
problems, in section \ref{subsec:recursive}, we characterize the
optimal switching times as exit times from certain intervals through
sequential optimal stopping problems equivalent to the original
switching problem.  In section \ref{subsec:value-function}, we shall
provide a new characterization of the value function, which leads to
a direct solution method described in \ref{subsec:method}.  We shall
illustrate this method through examples in section
\ref{sec:example}, one of which is a new optimal switching problem.
Section \ref{sec:last-section} concludes with comments on an
extension to a further general problem.

\section{Optimal Switching Problems}
We consider the following optimal switching problems for one
dimensional diffusions.  Let $(\Omega, \F, \p)$ be a complete
probability space with a standard Brownian motion $W=\{W_t; t\geq
0\}$.  Let $Z_t$ be the indicator vector at time $t$, $Z_t\in \{z_1,
z_2,..., z_m\}\triangleq \mathcal{Z}$ where each vector $z_i=(a_1,
a_2,..., a_k)$ with $a$ is either $0$ (closed) or $1$ (open), so
that $m=2^k$.  In this section, we consider the case of $k=1$. That
is, $Z_t$ takes either $0$ or $1$.  The admissible switching
strategy is
\begin{equation*}
w=(\theta_0, \theta_1, \theta_2,..., \theta_k,...; \zeta_0, \zeta_1,
\zeta_2,..., \zeta_k,...)
\end{equation*}
with $\theta_0=0$ where where where $0\leq \theta_1<\theta_2<....$
are an increasing sequence of $\mathcal{F}_t$-stopping times and
$\zeta_1$, $\zeta_2...$ are $\mathcal{F}_{\theta_i}$-measurable
random variables representing the new value of $Z_t$ at the
corresponding switching times $\theta_i$ (in this section,
$\zeta_i=1$ or $0$).  The state process at time $t$ is denoted by
$(X_t)_{t\ge 0}$ with state space $\mathcal{I}=(c, d)\subseteq \R$
and $X_0=x \in \mathcal{I}$,  and with the following dynamics:

If $\zeta_0=1$ (starting in open state), we have, for $m=0, 1,
2,.....$,
\begin{align}\label{eq:sde-open}
dX_t=\begin{cases} dX_t^1=\mu_1(X^1)dt + \sigma_1(X^1)dW_t,
&\theta_{2m}\le t<\theta_{2m+1},\\
dX_t^0=\mu_0(X^0)dt + \sigma_0(X^0)dW_t, &\theta_{2m+1}\le t <
\theta_{2m+2},
\end{cases}
\end{align}
and if $\zeta_0=0$ (starting in closed state),
\begin{align} \label{eq:sde-closed}
dX_t=\begin{cases} dX_t^0=\mu_0(X^1)dt + \sigma_0(X^0)dW_t,
&\theta_{2m}\le t<\theta_{2m+1},\\
dX_t^1=\mu_1(X^1)dt + \sigma_1(X^1)dW_t, &\theta_{2m+1}\le t <
\theta_{2m+2}.
\end{cases}
\end{align}
We assume that $\mu_i:\R\rightarrow \R$ and $\sigma_i:\R\rightarrow
\R$ are some Borel functions that ensure the existence and
uniqueness of the solution of (\ref{eq:sde-open}) for $i=1$ and
(\ref{eq:sde-closed}) for $i=0$.

Our performance measure, corresponding to starting state $i=0, 1$,
is
\begin{equation}\label{eq:problem}
J^w_i(x)=\ME\left[\int_0^\infty e^{-\alpha
s}f(X_s)ds-\sum_{j=1}^\infty e^{-\alpha \theta_j}H(X_{\theta_{j-}},
\zeta_j)\right]
\end{equation}
where $H: \mathbb{R}\times \mathcal{Z}\rightarrow \mathbb{R}_+$ is
the switching cost function and $f: \R\rightarrow \R$ is a
continuous function that satisfies
\begin{equation} \label{eq:f-condition}
\ME\left[\int_0^\infty e^{-\alpha s}|f(X_s)|ds\right]<\infty.
\end{equation}
In this section, the cost functions are of the form:
\begin{align*}
H(X_{\theta-}, \zeta)=  \begin{cases}
              H(X_{\theta-}, 1) &\text{opening cost},\\
              H(X_{\theta-}, 0) &\text{closing cost}.
              \end{cases}
\end{align*}
The optimal switching problem is to optimize the performance measure
for $i=0$ (start in closed state) and  $1$ (start in open state).
That is to find, for both $i=1$ and $i=0$,
\begin{eqnarray}\label{eq:problem-2}
v_i(x)\triangleq\sup_{w\in W}J^w(x) \quad\text{with}\quad X_0=x
\end{eqnarray}
where $W$ is the set of all the admissible strategies.
\subsection{Characterization of switching
times}\label{subsec:recursive} For the remaining part of section
$2$, we assume that the state space $X$ is $\mathcal{I}=(c, d)$
where both $c$ and $d$ are natural boundaries of $X$.  But our
characterization of the value function does not rely on this
assumption.  In fact, it is easily applied to other types of
boundaries, for example, absorbing boundary.

The first task is to characterize the optimal switching times as
exit times from intervals in $\R$. For this purpose, we define two
functions $g_0$ and  $g_1: \R_+\rightarrow \R$ with
\begin{equation}\label{eq:g-1-g-0}
g_1(x)\triangleq \sup_{w\in W_0}J^w_1(x)
\quad\text{and}\quad
g_0(x)\triangleq\sup_{w\in W_0}J^w_0(x).
\end{equation}
where $W_0\triangleq \{w\in W: w=(\theta_0, \zeta_0,
\theta_1=+\infty )\}$. In other words, $g_1(\cdot)$ is the
discounted expected revenue by starting with $\zeta_0=1$ and making
no switches. Similarly, $g_0(\cdot)$ is the discounted expected
revenue by staring with $\zeta_0=0$ and making no switches.

We set $w_0\triangleq g_1$ and $y_0\triangleq g_0$. We consider the
following simultaneous sequential optimal stopping problems with
$w_n : \R_+ \rightarrow \R$ and $y_n: \R_+\rightarrow \R$ for $n=1,
2,....$:
\begin{equation}\label{eq:w-function}
w_{n}(x)\triangleq\sup_{\tau\in \mathcal{S}}\ME \left[\int_0^\tau
e^{-\alpha s}f(X_s)ds + e^{-\alpha\tau}
(y_{n-1}(X_{\tau})-H(X_{\tau-}, 1-Z_{\tau-}))\right],
\end{equation}
and
\begin{equation}\label{eq:y-function}
y_{n}(x)\triangleq\sup_{\tau\in \mathcal{S}}\ME\left[\int_0^\tau
e^{-\alpha s}f(X_s)ds + e^{-\alpha\tau}
(w_{n-1}(X_{\tau})-H(X_{\tau-}, 1-Z_{\tau-}))\right],
\end{equation}
where $\mathcal{S}$ is a set of $\mathcal{F}_t$ stopping times. Note
that for each $n$, the sequential problem \ref{eq:w-function} (resp.
(\ref{eq:y-function})) starts in open (resp. closed) state.

On the other hand, we define $n$-time switching problems for
$\zeta_0=1$:
\begin{equation}\label{eq:q-function}
q^{(n)}(x)\triangleq\sup_{w\in W_n}J_1^w(x),
\end{equation}
where \begin{equation*}
    W_n\triangleq \{w\in W; w=(\theta_1, \theta_2,... \theta_{n+1}; \zeta_1,
    \zeta_2,...\zeta_n); \theta_{n+1}=+\infty\}.
\end{equation*}
In other words, we start with $\zeta_0=1$ (open) and are allowed to
make at most $n$ switches.  Similarly, we define another $n$-time
switching problems corresponding to $\zeta_0=0$:
\begin{equation}\label{eq:z-function}
p^{(n)}(x)\triangleq\sup_{w\in W_n}J_0^w(x).
\end{equation}
We investigate the relationship of these four problems:
\begin{lemma}
For any $x \in \R$, $w_{n}(x)=q^{(n)}(x)$ and $y_{n}(x)=p^{(n)}(x)$.
\end{lemma}
\begin{proof}
We shall prove only the first assertion since the proof of the
second is similar.  We have set $y_0(x)=g_0(x)$. Now we consider
$w_1$ by using the strong Markov property of $X$:
\begin{align*}
w_1(x)&=\sup_{\tau\in \S}\ME\left[\int_0^\tau e^{-\alpha s}f(X_s)ds
+ e^{-\alpha\tau} (g_{0}(X_{\tau})-H(X_{\tau-}, 0))\right]\\
&=\sup_{\tau\in \S}\ME\left[\int_0^\infty e^{-\alpha s}f(X_s)ds
-\int_{\tau}^\infty e^{-\alpha s}f(X_s)ds-e^{-\alpha
\tau}(g_0(X_{\tau})-H(X_{\tau-}, 0))
\right]\\
&=\sup_{\tau\in \S}\ME\left[e^{-\alpha
\tau}(g_0(X_{\tau})-g_1(X_{\tau})-H(X_{\tau-}, 0))\right]+g_1(x).
\end{align*}
On the other hand,
\begin{align*}
q^{(1)}(x)&=\sup_{w\in W_1}\ME\left[\int_{0}^{\infty}e^{-\alpha
s}f(X_s)ds -e^{-\alpha \theta_1}H(X_{\theta_{1-}},
\zeta_1)\right]\\
&=\sup_{w\in W_1}\ME\left[\int_{0}^{\theta_1}e^{-\alpha
s}f(X_s)ds+\int_{\theta_1}^{\infty}e^{-\alpha s}f(X_s)ds -e^{-\alpha
\theta_1}H(X_{\theta_{1-}}, 0)\right]\\
&=\sup_{w\in
W_1}\ME\left[(g_1(x)-e^{-\alpha\theta_1}g_1(X_{\theta_1}))-e^{-\alpha
\theta_1}(g_0(X_{\theta_1})-H(X_{\theta_{1-}}, 0))\right]\\
&=\sup_{w\in
W_1}\ME\left[e^{-\alpha\theta_1}(g_0(X_{\theta_1})-g_1(X_{\theta_1})-H(X_{\theta_{1-}},
0))\right]+g_1(x).
\end{align*}
Since both $\tau$ and $\theta_1$ are $\F_t$ stopping times, we have
$w_1(x)=q^{(1)}(x)$ for all $x\in \R$.  Moreover, by the theory of
the optimal stopping (see Appendix \ref{appx:ost}, especially
Proposition \ref{prop:A4}), $\tau$ and hence $\theta_1$ are
characterized as an exit time from an interval.  Similarly, we can
prove $y_1(x)=p^{(1)}(x)$. Now we consider $q^{(2)}(x)$ which is the
value if we start in open state and make at most $2$  switches (open
$\rightarrow$ close $\rightarrow$ open).
For this purpose, we consider the performance measure
$\bar{q}^{(2)}$ that starts in an open state and is allowed two
switches:  For arbitrary switching times $\theta_1,
\theta_2>\theta_1 \in \S$, we have
\begin{align*}
\bar{q}^{(2)}(x)&\triangleq\ME\left[\int_{0}^{\infty}e^{-\alpha
s}f(X_s)ds -\sum_{j=1}^2e^{-\alpha \theta_j}H(X_{\theta_{j-}},
\zeta_j)\right]\\
&=\ME\Bigg[\int_{0}^{\theta_1}e^{-\alpha s}f(X_s)ds+\int_{\theta_1}^{\theta_2}e^{-\alpha s} f(X_s)ds+\int_{\theta_2}^{\infty}e^{-\alpha s} f(X_s)ds\\
&\hspace{3cm}-e^{-\alpha\theta_1}H(X_{\theta_1-},
0)-e^{-\alpha\theta_2}H(X_{\theta_2-}, 1)\Bigg]\\
&=\left(g_1(x)-\ME[e^{-\alpha\theta_1}g_1(X_{\theta_1})]\right)
+\left(\ME[e^{-\alpha\theta_1}g_0(X_{\theta_1})-e^{-\alpha\theta_2}g_0(X_{\theta_2})]\right)
+\ME[e^{-\alpha\theta_2}g_1(X_{\theta_2})]\\
&\hspace{3cm}-\ME[e^{-\alpha\theta_1}H(X_{\theta_1-},
0)+e^{-\alpha\theta_2}H(X_{\theta_2-}, 1)].
\end{align*}
Hence we have the following multiple optimal stopping problems:
\begin{align*}
\bar{q}^{(2)}(x)&=\sup_{(\theta_1, \theta_2)\in
\S^2}\ME\left[e^{-\alpha\theta_1}\Big((g_0-g_1)(X_{\theta_1})-H(X_{\theta_1-},
0)\Big)+e^{-\alpha\theta_2}\Big((g_1-g_0)(X_{\theta_2})-H(X_{\theta_2-},
1)\Big)\right]\\
&\hspace{3cm}+g_1(x)
\end{align*}
where $\S^2\triangleq \{(\theta_1, \theta_2); \theta_1\in\S;
\theta_2\in\S_{\theta_1}\}$ and $\S_{\sigma}=\{\tau\in\S; \tau\ge
\sigma\}$ for every $\sigma\in\S$.  Let us denote $h_1(x)\triangleq
g_1(x)-g_0(x)-H(x, 0)$,   
$h_2(x)\triangleq g_0(x)-g_1(x)-H(x, 1)$,
\begin{align*}
V_1(x)\triangleq\sup_{\tau\in\S}\ME\left[e^{-\alpha\tau}h_1(X_{\tau})\right]\quad\text{and}\quad
V_2(x)\triangleq\sup_{\tau\in\S}\ME\left[e^{-\alpha\tau}(h_2(X_\tau)+V_1(X_{\tau}))\right].
\end{align*}
We also define
\begin{equation*}
  \Gamma_1\triangleq \{x\in \mathcal{I}:
  V_1(x)=h_1(x)\}\quad\text{and}\quad\Gamma_2\triangleq
  \{x\in\mathcal{I}: V_2(x)=h_2(x)+V_1(x)\}
\end{equation*}
with $\sigma_n\triangleq \inf\{t\ge 0: X_t\in\Gamma_n\}$.  By using
Proposition 5.4. in Carmona and Dayanik \cite{CD2003},  we
conclude that $\theta_1=\sigma_1$ and
$\theta_2=\theta_1+\sigma_2\circ s(\theta_1)$ is optimal strategy
where $s(\cdot)$ is the shift operator.  Hence we only consider the
maximization over the set of admissible strategy $W_2^*$ where
\begin{equation*}
W^{*}_2\triangleq\{w\in W_2:  \theta_1, \theta_2 \quad\text{are exit
imes from an interval in $\mathcal{I}$}\},
\end{equation*}
and can use the relation $\theta_2-\theta_1=\theta \circ
s(\theta_1)$ with some exit time $\theta\in\S$.
\begin{align*}
q^{(2)}(x)&=\sup_{w\in W_2^*}\ME\left[\int_{0}^{\infty}e^{-\alpha
s}f(X_s)ds -\sum_{j=1}e^{-\alpha \theta_j}H(X_{\theta_{j-}},
\zeta_j)\right]\\
&=\sup_{w\in W_2^*}\ME\Big[\int_{0}^{\theta_1}e^{-\alpha s}f(X_s)ds+\int_{\theta_1}^{\theta_2}e^{-\alpha s} f(X_s)ds+\int_{\theta_2}^{\infty}e^{-\alpha s} f(X_s)ds\\
&\hspace{3cm}-e^{-\alpha\theta_1}(H(X_{\theta_{1-}},
0)+e^{-\alpha(\theta_2-\theta_1)} H(X_{\theta_{2-}}, 1))\Big]\\
&=\sup_{w\in W_2^*}\ME\Big[\int_{0}^{\theta_1}e^{-\alpha
s}f(X_s)ds+e^{-\alpha\theta_1}\E^{X_{\theta_1}}\left[\left(\int_0^{\theta}
+\int_\theta^{\infty}\right)e^{-\alpha
s}f(X_s)ds-e^{-\alpha\theta}H(X_{\theta-}, 1)\right]\\
&\hspace{3cm}-e^{-\alpha\theta_1}H(X_{\theta_{1-}}, 0)\Big].
\end{align*}
Now by using the result for $p^{(1)}$, we can conclude
\begin{align*}
q^{(2)}(x)&=\sup_{w\in W_2^*}\ME\left[\int_{0}^{\theta_1}e^{-\alpha
s}f(X_s)ds+e^{-\alpha\theta_1}\left(p^{(1)}(X_{\theta_1})-H(X_{\theta_{1-}}, 0)\right)\right]\\
&=\sup_{\theta_1\in \S}\ME\left[\int_{0}^{\theta_1}e^{-\alpha
s}f(X_s)ds+e^{-\alpha\theta_1}\left(y_1(X_{\theta_1})-H(
X_{\theta_{1-}}, 0)\right)\right] =w_2(x)
\end{align*}
Similarly, we can prove $y_2(x)=p^{(2)}(x)$ and we can continue this
process inductively to conclude that $w_{n}(x)=q^{(n)}(x)$ and
$y_{n}(x)=p^{(n)}(x)$ for all $x$ and $n$.
\end{proof}
\begin{lemma} For all $x \in \R$,
$\lim_{n\rightarrow \infty}q^{(n)}(x)=v_1(x)$ and
$\lim_{n\rightarrow \infty}p^{(n)}(x)=v_0(x)$.
\end{lemma}
\begin{proof}
Let us define $q(x)\triangleq \lim_{n\rightarrow \infty}q^{(n)}(x)$.
Since $W_n\subset W$, $q^{(n)}(x)\le v_1(x)$ and hence $q(x)\le
v_1(x)$. To show the reverse inequality, we define $W^+$ to be a set
of admissible strategies such that
\begin{equation*}
  W^+=\{w\in W: J_1^w(x)< \infty \quad\text{for all}\quad x\in \R\}.
\end{equation*}
Let us assume that $v_1(x)<+\infty$ and consider a strategy $w^+\in
W^+$ and another strategy $w_n$ that coincides with $w^+$ up to and
including time $\theta_n$ and then takes no further interventions.
\begin{equation}\label{eq:difference}
  J^{w^+}_1(x)-J^w_1(x)=\ME\left[\int_{\theta_n}^\infty e^{-\alpha
  s}(f(X_s)-f(X_{s-{\theta_n}}))
  -\sum_{i\ge n+1}e^{-\alpha \theta_i}H(X_{\theta_i-},
  \zeta_i)\right],
\end{equation}
which implies
\begin{equation*}
  |J^{w^+}_1(x)-J^w_1(x)|\le \ME\left[\frac{2\|f\|}{\alpha}e^{-\alpha\theta_n}-\sum_{i\ge n+1}e^{-\alpha \theta_i}H(X_{\theta_i-},
  \zeta_i)\right].
\end{equation*}
As $n\rightarrow +\infty$, the right hand side goes to zero by the
dominated convergence theorem.  Hence it is shown
\begin{equation*}
  v_1(x)=\sup_{w\in W^+}J^w_1(x)=\sup_{w\in {\cup}_n W_n}J^w_1(x)
\end{equation*}
so that $v_1(x)\le q(x)$.  Next we consider $v_1(x)=+\infty$.  Then
we have some $m\in \mathbb{N}$ such that $w_m(x)=q^{(m)}(x)=\infty$.
Hence $q^{(n)}(x)=\infty$ for all $n\ge m$.  The second assertion is
proved similarly.
\end{proof}
We define an operator $\LL: \mathcal{H}\rightarrow \mathcal{H}$
where $\mathcal{H}$ is a set of Borel functions
\begin{align*}
\LL u(x)\triangleq \sup_{\tau\in\S}\ME\left[\int_0^\tau e^{-\alpha
s}f(X_s)ds +e^{-\alpha\tau}\left(u(X_{\tau})-H(X_{\tau-},
1-Z_{\tau-})\right)\right].
\end{align*}
\begin{lemma}\label{lem:fixed-point}
The function $w(x)\triangleq\lim_{n\rightarrow \infty}w_n(x)$ is the
smallest solution, that majorizes $g_1(x)$, of the function equation
$w=\LL w$.
\end{lemma}
\begin{proof}
We renumber the sequence $(w_0, y_1, w_2, y_3...)$ as $(u_0, u_1,
u_2, u_3....)$.  Since $u_n$ is monotone increasing, the limit
$u(x)$ exists. We have $u_{n+1}(x)=\LL u_{n}(x)$ and apply the
monotone convergence theorem by taking $n\rightarrow \infty$, we
have $u(x)=\LL u(x)$. We assume that $u'(x)$ satisfies $u'=\LL u'$
and majorizes $g_1(x)=u_0(x)$. Then $u'=\LL u'\geq \LL u_0=u_1$. Let
us assume, for induction argument that $u'\geq u_n$, then
\begin{align*}
  u'=\LL u'\geq \LL u_n=u_{n+1}.
\end{align*}
Hence we have $u'\geq u_n$ for all $n$, leading to $u'\geq
\lim_{n\rightarrow \infty}u_n=u$.  Now we take the subsequence in
$(w_0, y_1, w_2, y_3....)$ to complete the proof.
\end{proof}

\begin{proposition}
For each $x\in \R$, $\lim_{n\rightarrow \infty}w_n(x)=v_1(x)$ and
$\lim_{n\rightarrow \infty}y_n(x)=v_0(x)$.  Moreover, the optimal
switching times, $\theta^*_i$ are exit times from an interval.
\end{proposition}
\begin{proof}
 We can prove the first
assertion by combining the first two lemmas above. Now we
concentrate on the sequence of $w_n(x)$. For each $n$, finding
$w_{n}(x)$ by solving (\ref{eq:w-function}) is an optimal stopping
problem.  By Proposition \ref{prop:A4}, the optimal stopping times
are characterized as an exit time of $X$ from an interval for all
$n$. This is also true in the limit: Indeed, by Lemma
\ref{lem:fixed-point}, in the limit, the value function of optimal
switching problem $v_1(x)=w(x)$ satisfies $w=\LL w$, implying that
$v_1(x)$ is the solution of an optimal stopping problem. Hence the
optimal switching times are characterized as exit time from an
interval.
\end{proof}
\subsection{Characterization of the value
functions}\label{subsec:value-function} We go back to the original
problem (\ref{eq:problem}) to characterize the value function of the
optimal switching problems.
By the exit time characterization of the optimal switching times,
$\theta_i^*$ are given by
\begin{align}\label{eq:optimal-switching-time} \theta_i^*=\begin{cases}\inf\{t>\theta_{i-1}; X_t^1\in
\Gamma_1\}\\
\inf\{t>\theta_{i-1}; X_t^0\in \Gamma_0\} \end{cases}
\end{align}
where $\Gamma_1=\R\setminus \mathbf{C}_1$ and $\Gamma_0=\R\setminus
\mathbf{C}_0$.  We define here $\mathbf{C}_i$ and $\Gamma_i$ to be continuation and stopping region for $X_t^i$, respectively.  We can simplify the performance measure
$J^w$ considerably. For $\zeta_0=1$, we have
\begin{align*}
J^w_1(x)&=\ME\left[\int_0^\infty e^{-\alpha
s}f(X_s)ds-\sum_{j=1}^\infty e^{-\alpha
\theta_j}H(X_{\theta_{j-}}, \zeta_j)\right]\\
&=\ME\Bigg[\int_0^{\theta_1}e^{-\alpha s}f(X_s)ds +
\int_{\theta_1}^\infty e^{-\alpha
s}f(X_s)ds\\
&\hspace{3cm}-e^{-\alpha\theta_1}\left(H(X_{\theta_1-},
0)+\sum_{j=2}e^{-\alpha(\theta_i-\theta_1)}H(X_{\theta_{j-}},
\zeta_j)\right) \Bigg]\\
&=\ME\Bigg[\int_0^{\theta_1}e^{-\alpha
s}f(X_s)ds+e^{-\alpha\theta_1}\E^{X_{\theta_1}}\left[\int_0^{\infty}e^{-\alpha
s}f(X_s)ds -\sum_{j=1}e^{-\alpha\theta_j}H(X_{\theta_{j-}},
\zeta_j)\right]\\
&\hspace{3cm}-e^{-\alpha\theta_1}H(X_{\theta_1-}, 0)\Bigg]
\end{align*}
We notice that in the time interval $(0, \theta_1)$, the process $X$
is not intervened. The inner expectation is just
$J_0^w(X_{\theta_1})$.  Hence we further simplify
\begin{align*}
J_1^w(x)&=\ME\left[\int_0^{\theta_1}e^{-\alpha
s}f(X_s)ds+e^{-\alpha\theta_1}(J_0^w(X_{\theta_1})-H(X_{\theta_1-},
0))\right]\\
&=\ME\left[-e^{-\alpha\theta_1}g_1(X_{\theta_1})+e^{-\alpha\theta_1}(J_0^w(X_{\theta_1})-H(X_{\theta_1-},
0))\right]+g_1(x)\\
&=\ME\left[-e^{-\alpha\theta_1}g_1(X_{\theta_1})+e^{-\alpha\theta_1}J_1^w(X_{\theta_1})\right]+g_1(x).
\end{align*}
The third equality is a critical observation.  Finally, we define
$u_1\triangleq J_1-g_1$ and obtain
\begin{equation}\label{eq:u-1}
u_1(x)=J^w_1(x)-g_1(x)=\ME\left[e^{-\alpha\theta_1}u_1(X_{\theta_1})\right].
\end{equation}
Since the switching time $\theta_1$ is characterized as a hitting
time of a certain point in the state space, we can represent
$\theta_1=\tau_a\triangleq\inf\{t\geq 0: X_t=a\}$ for some $a\in
\R$. Hence equation (\ref{eq:u-1}) is an optimal stopping problem
that maximizes
\begin{equation}
u_1(x)=J^w_1(x)-g_1(x)=\ME\left[e^{-\alpha\tau_a}u_1(X_{\tau_a})\right].
\end{equation}
among all the $\tau_a\in \S$.  When $\theta_1=0$ (i.e.,
$x=X_{\theta_1}$),
\begin{align*}
J_1^w(x)&=\ME\left[-g_1(x)+J_0^w(x)-H(x, 0)\right]+g_1(x)
\end{align*}
and hence
\begin{align*}
u_1(x)&=J_0^w(x)-H(x, 0)-g_1(x).
\end{align*}
In other words, we make a switch from open to closed immediately  by
paying the switching cost. Similarly, for $\zeta_0=0$, we can
simplify the performance measure $J_0^w(\cdot)$ to obtain
\begin{equation*}
J_0^w(x)=\ME\left[-e^{-\alpha\theta_1}g_0(X_{\theta_1})+e^{-\alpha\theta_1}J_0^w(X_{\theta_1})\right]+g_0(x).
\end{equation*}
By defining $u_0\triangleq J_0^w-g_0$, we have
\begin{equation*}
u_0(x)=J^w_0(x)-g_0(x)=\ME\left[e^{-\alpha\theta_1}u_0(X_{\theta_1})\right].
\end{equation*}
Again, by using the characterization of switching times, we replace
$\theta_1$ with $\tau_b$,
\begin{equation}\label{eq:u-0}
u_0(x)=J^w_0(x)-g_0(x)=\ME\left[e^{-\alpha\tau_b}u_0(X_{\tau_b})\right].
\end{equation}
In summary, we have
\begin{align}\label{eq:ost-1}
u_1(x)&=\begin{cases} u_0(x)+g_0(x)-H(x, 0)-g_1(x), &x\in
\Gamma_1,\\
\ME\left[e^{-\alpha\tau_a}u_1(X_{\tau_a})\right]=\ME\left[e^{-\alpha\tau_a}(u_0(X_{\tau_a})+g_0(X_{\tau_a})-
g_1(X_{\tau_a})-H(X_{\tau_a}, 0))\right], &x\in \mathbf{C}_1,
\end{cases}
\end{align}
and
\begin{align}
u_0(x)&=\begin{cases}\label{eq:ost-0}
\ME\left[e^{-\alpha\tau_b}u_0(X_{\tau_b})\right]=\ME\left[e^{-\alpha\tau_b}(u_1(X_{\tau_b})+g_1(X_{\tau_b})-g_0(X_{\tau_b})-H(X_{\tau_b}, 1))\right], &x\in \mathbf{C}_0,\\
u_1(x)+g_1(x)-H(x, 1)-g_0(x), &x\in \Gamma_0.
\end{cases}
\end{align}
Hence we should solve the following optimal stopping problems
simultaneously:
\begin{align}\label{eq:system}
\begin{cases}
\bar{v}_1(x)\triangleq\sup_{\tau\in\S}\ME\left[e^{-\alpha\tau}(u_1(X_{\tau})\right]\\
\bar{v}_0(x)\triangleq
\sup_{\sigma\in\S}\ME\left[e^{-\alpha\sigma}(u_0(X_{\sigma})\right]
\end{cases}
\end{align}

Now we let the infinitesimal generators of $X^1$ and $X^0$ be $\A_1$
and $A_0$, respectively.   We consider $(\A_i-\alpha)v(x)=0$ for
$i=0, 1$. This ODE has two fundamental solutions, $\psi_i(\cdot)$
and $\varphi_i(\cdot)$. We set $\psi_i(\cdot)$ is an increasing and
$\varphi_i(\cdot)$ is a decreasing function. Note that
$\psi_i(c+)=0, \varphi_i(c+)=\infty$ and $\psi_i(d-)=\infty,
\varphi_i(d-)=0$. We define
\begin{equation*}
F_i(x)\triangleq\frac{\psi_i(x)}{\varphi_i(x)}\quad \text{and}\quad
G_i(x)\triangleq -\frac{\varphi_i(x)}{\psi_i(x)}\quad \text{for
$i=0, 1$}.
\end{equation*}
\noindent By referring to Dayanik and Karatzas~\cite{DK2003}, we
have the following representation
\begin{equation*}
\ME[e^{-\alpha\tau_r}1_{\{\tau_r<\tau_l\}}]=\frac{\psi(l)\varphi(x)-\psi(x)\varphi(l)}
{\psi(l)\varphi(r)-\psi(r)\varphi(l)},\quad
\ME[e^{-\alpha\tau_r}1_{\{\tau_l<\tau_r\}}]=\frac{\psi(x)\varphi(r)-\psi(r)\varphi(x)}
{\psi(l)\varphi(r)-\psi(r)\varphi(l)},
\end{equation*}
for $x\in[l,r]$ where $\tau_l\triangleq\inf\{t>0; X_t=l\}$ and
$\tau_r\triangleq\inf\{t>0; X_t=r\}$.

By defining
\begin{equation*}
 W_1=(u_1/\psi_1)\circ G^{-1}_1 \quad \text{and}\quad
W_0=(u_0/\varphi_0)\circ F^{-1}_0,
\end{equation*}
 the second equation in
(\ref{eq:ost-1}) and the first equation in (\ref{eq:ost-0}) become
\begin{align}\label{eq:W1}
W_1(G_1(x))&=W_1(G_1(a))\frac{G_1(d)-G_1(x)}{G_1(d)-G_1(a)}+
W_1(G_1(d))\frac{G_1(x)-G_1(a)}{G_1(d)-G_1(a)} \quad x\in[a, d),
\end{align}
and
\begin{align}\label{eq:W0}
W_0(F_0(x))&=W_0(F_0(c))\frac{F_0(b)-F_0(x)}{F_0(b)-F_0(c)}+W_0(F_0(b))\frac{F_0(x)-F_0(c)}{F_0(b)-F_0(c)},
\quad x\in(c, b],
\end{align}
respectively.  We should understand that $F_0(c)\triangleq
F_0(c+)=\psi_0(c+)/\varphi_0(c+)=0$ and that $G_1(d)\triangleq
G_1(d-)=-\varphi_1(d-)/\psi_1(d-)=0$. In the next subsection, we
shall explain $W_1(G_1(d-))$ and $W_0(F_0(c+))$ in details.  Both
$W_1$ and $W_0$ are a \emph{linear function} in their respective
transformed spaces. Hence under the appropriate transformations, the
two value functions are linear functions in the continuation region.

\subsection{Direct Method for a Solution}\label{subsec:method}
We have established a mathematical characterization of the value
functions of optimal switching problems.  We shall investigate, by
using the characterization, a direct solution method that does not
require the recursive optimal stopping schemes described in section
\ref{subsec:recursive}.  Since the two optimal stopping problems
(\ref{eq:system}) have to be solved simultaneously, finding $u_0$ in
$x\in \C_0$, for example, requires that we find the smallest
$F_0$-concave majorant of $(u_1(x)+g_1(x)-g_0(x)-H(x,
1))/\varphi_0(x)$ as in (\ref{eq:ost-0}) that involves $u_1$.

There are two cases, depending on whether $x \in \mathbf{C}_1\cap
\C_0$ or $x\in \Gamma_1\cap \C_0$, as to what $u_1(\cdot)$
represents.  In the region $x\in \Gamma_1 \cap \C_0$, $u_1(\cdot)$
that shows up in the equation of $u_0(x)$ is of the form
$u_1(x)=u_0(x)+g_0(x)-H(x, 1, 0)-g_1(x)$. In this case, the
``obstacle" that should be majorized is in the form
\begin{align}\label{eq:on-(c, a)}
&u_1(x)+g_1(x)-g_0(x)-H(x, 1)\nonumber\\
&=(u_0(x)+g_0(x)-H(x, 0)-g_1(x))+g_1(x)-g_0(x)-H(x, 1)\nonumber\\
&=u_0(x)-H(x, 0)-H(x, 1)<u_0(x).
\end{align}
This implies that in $x\in \Gamma_1 \cap \C_0$, the $u_0(x)$
function always majorizes the obstacle.  Similarly, in $x\in
\Gamma_0\cap \C_1$, the $u_1(x)$ function always majorizes the
obstacle.

Next, we consider the region $x\in \mathbf{C}_0 \cap \C_1$.  
The $u_0(\cdot)$ term in (\ref{eq:ost-1}) is represented, due to its
linear characterization, as
\begin{equation*}
W_0(F_0(x))=\beta_0(F_0(x))+ d_0
\end{equation*}
with some $\beta_0\in \R$ and  $d_0 \in \R_+$ in the transformed
space. (The nonnegativity of $d_0$ will be shown.) In the original
space, it has the form of $\varphi_0(x)(\beta_0F_0(x)+d_0)$. Hence
by the transformation $(u_1/\psi_1)\circ G^{-1}$, $W_1(G_1(x))$ is
the smallest linear majorant of
\begin{equation*}
\frac{K_1(x)+\varphi_0(x)(\beta_0F_0(x)+d_0)}{\psi_1(x)}=\frac{K_1(x)+\beta_0\psi_0(x)+d_0\varphi_0(x)}{\psi_1(x)}
\end{equation*}
on $(G_1(d-), G_1(a^*))$ where
\begin{equation}
K_1(x)\triangleq
g_0(x)-g_1(x)-H(x, 0).
\end{equation}
 This linear function passes a point
$(G_1(d-), l_d)$ where $G_1(d-)=0$ and
\begin{equation*}
l_d=\limsup_{x\uparrow
d}\frac{(K_1(x)+\beta_0\psi_0(x)+d_0\varphi_0(x))^+}{\psi_1(x)}.
\end{equation*}
Let us consider further the quantity $l_d \geq 0$.  By noting
\begin{align*}
\limsup_{x\uparrow d}\frac{(K_1(x)+\beta_0\psi_0(x))^+}{\psi_1(x)}
&\leq \limsup_{x\uparrow
d}\frac{(K_1(x)+\beta_0\psi_0(x)+d_0\varphi_0(x))^+}{\psi_1(x)}\\
&\leq \limsup_{x\uparrow
d}\frac{(K_1(x)+\beta_0\psi_0(x))^+}{\psi_1(x)}+\limsup_{x\uparrow
d}\frac{d_0\varphi_0(x)}{\psi_1(x)}
\end{align*}
and $\limsup_{x\uparrow d}\frac{d_0\varphi_0(x)}{\psi_1(x)}=0$, we
can redefine $l_d$ by
\begin{align}\label{eq:ld}
l_d\triangleq \limsup_{x\uparrow
d}\frac{(K_1(x)+\beta_0\psi_0(x))^+}{\psi_1(x)}
\end{align}
to determine the finiteness of the value function of the optimal
switching problem, $v_1(x)$, based upon Proposition
\ref{prop:A5}-\ref{prop:A7}. 
Let us  concentrate on the case $l_d=0$.

Similar analysis applies to (\ref{eq:ost-0}). $u_1(x)$ in
(\ref{eq:ost-0}) is represented as
\begin{equation*}
W_1(G_1(x))=\beta_1G_1(x)+d_1
\end{equation*}
with some $\beta_1\in \R$ and $d_1\in\R_+$.  Note that $d_1=l_d\geq 0$.  
In the original space, it has the form of
$\psi_1(x)(\beta_1G_1(x)+d_1)$.  Hence by the transformation
$(u_0/\varphi_0(x))\circ F^{-1}$, $W_0(F_0(x))$ is the smallest
linear majorant of
\begin{equation*}
\frac{K_0(x)+\psi_1(x)(\beta_1G_1(x)+d_1)}{\varphi_0(x)}=\frac{K_0(x)-\beta_1\varphi_1(x)+d_1\psi_1(x)}{\varphi_0(x)}
\end{equation*}
on $(F_0(c+), F_0(b^*))$ where
\begin{equation}
K_0(x)\triangleq g_1(x)-g_0(x)-H(x, 1).
\end{equation}  This linear
function passes a point $(F_0(c+), l_c)$ where $F_0(c+)=0$ and
\begin{equation*}
l_c=\limsup_{x\downarrow
c}\frac{(K_0(x)-\beta_1\varphi_1(x)+d_1\psi_1(x))^+}{\varphi_0(x)}.
\end{equation*}
Hence we have $l_c=d_0\geq 0$.  By the same argument as for $l_d$,
we can redefine
\begin{align}\label{eq:lc}
l_c\triangleq \limsup_{x\downarrow
c}\frac{(K_0(x)-\beta_1\varphi_1(x))^+}{\varphi_0(x)}.
\end{align}

\begin{remark}\normalfont
\begin{itemize}\label{rem:lc-ld}
\item[(a)]
Evaluation of $l_d$ or $l_c$ does not require knowledge of $\beta_0$
or $\beta_1$, respectively unless the orders of $\max(K_1(x),
\psi_1(x))$ and $\psi_0(x)$ are equal, for example. (For this event,
see Proposition \ref{prop:odd-case}.) Otherwise, we just compare the
order of the positive leading terms of the numerator in
(\ref{eq:ld}) and (\ref{eq:lc}) with that of the denominator.

\item[(b)] \emph{A sufficient condition for $l_d=l_c=0$}: since we have
\begin{equation*}
0\leq l_d\leq \limsup_{x\uparrow
d}\frac{(K_1(x))^+}{\psi_1(x)}+\limsup_{x\uparrow
d}\frac{(\beta_0\psi_0(x))^+}{\psi_1(x)}.
\end{equation*}
a sufficient condition for $l_d=0$ is
\begin{equation}\label{eq:sufficent-d}
\limsup_{x\uparrow d}\frac{(K_1(x))^+}{\psi_1(x)}=0 \quad
\text{and}\quad \limsup_{x\uparrow d}\frac{\psi_0(x)}{\psi_1(x)}=0.
\end{equation}
Similarly,
\begin{equation*} 0\leq
l_c\leq \limsup_{x\downarrow
c}\frac{(K_0(x))^+}{\varphi_0(x)}+\limsup_{x\downarrow
c}\frac{(-\beta_1\varphi_1(x))^+}{\varphi_0(x)}
\end{equation*}
Hence a sufficient condition for $l_c=0$ is
\begin{equation}\label{eq:sufficient-c}
\limsup_{x\downarrow c}\frac{(K_0(x))^+}{\varphi_0(x)}=0\quad\text
{and}\quad \limsup_{x\downarrow
c}\frac{\varphi_1(x)}{\varphi_0(x)}=0.
\end{equation}
Moreover, it is obvious $\beta_1<0$ and $\beta_0>0$ since the linear
majorant passes the origin of each transformed space.  Recall a
points in the interval $(c, d)\in\R_+$ will be transformed by
$G(\cdot)$ to $(G(c), G(d-))\in \R_-$.
\end{itemize}
\end{remark}
\noindent We summarize the case of $l_c=l_d=0$:
\begin{proposition}\label{prop:2}
Suppose that $l_d=l_c=0$, the quantities being defined by
(\ref{eq:ld}) and by (\ref{eq:lc}), respectively.  The value
functions in the transformed space are the smallest linear majorants
of
\begin{equation*}
R_1(\cdot)\triangleq
\frac{r_1(G_1^{-1}(\cdot))}{\psi_1(G_1^{-1}(\cdot))}\quad\text{and}\quad
R_0(\cdot)\triangleq
\frac{r_0(F_0^{-1}(\cdot))}{\varphi_0(F_0^{-1}(\cdot))}
\end{equation*}
where
\begin{equation*}
r_1(x)\triangleq g_0(x)-g_1(x)+\beta_0\psi_0(x)-H(x, 0)
\end{equation*}
and
\begin{equation*}
r_0(x)\triangleq g_1(x)-g_0(x)-\beta_1\varphi_1(x)-H(x, 1)
\end{equation*}  for
\begin{equation}
\beta_0>0\quad\text{and}\quad \beta_1<0.
\end{equation}
 Furthermore, $\Gamma_1$ and $\Gamma_0$ in
(\ref{eq:ost-1}) and (\ref{eq:ost-0}) are given by
\begin{equation*}
\Gamma_1\triangleq \{x\in (c, d): W_1(G_1(x))=R_1(G_1(x))\},
\quad\text{and}\quad
\Gamma_0\triangleq \{x\in (c, d): W_0(F_0(x))=R_0(F_0(x))\}.
\end{equation*}
\end{proposition}

\begin{coro}\normalfont
If either of the boundary points $c$ or $d$ is \emph{absorbing},
then $(F_0(c), W_0(F_0(c))$ or $(G_1(d), W_1(G_1(d)))$ is obtained
directly.  We can entirely omit the analysis of $l_c$ or $l_d$. The
characterization of the value function (\ref{eq:W1}) and
(\ref{eq:W0}) remains exactly the same.
\end{coro}

\begin{remark}\normalfont
An algorithm to find $(a^*, b^*, \beta_0^*, \beta_1^*)$ can be
described as follows:
\begin{enumerate}
\item Start with some $\beta_1'\in \R$.
\item Calculate $r_0$ and then $R_0$ by the transformation
$R_0(\cdot)=
\frac{r_0(F_0^{-1}(\cdot))}{\varphi_0(F_0^{-1}(\cdot))}$.
\item Find the linear majorant of $R_0$ passing the origin of the
transformed space. Call the slope of the linear majorant, $\beta_0$
and the point, $F_0(b)$, where $R_0$ and the linear majorant meet .
\item Plug $b$ and $\beta_0$ in the equation for $r_1$ and calculate
$R_1$ by the transformation $R_1(\cdot)=
\frac{r_1(G_1^{-1}(\cdot))}{\psi_1(G_1^{-1}(\cdot))}$.
\item Find the linear majorant of $R_1$ passing the origin of the
transformed space.  Call the slope of the linear majorant, $\beta_1$
and the point, $G_1(a)$, where $R_1$ and the linear majorant meet.
\item Iterate step 1 to 5 until $\beta_1=\beta_1'$.
\end{enumerate}
If both $R_1$ and $R_0$ are differentiable functions with their
respective arguments, we can find $(a^*, b^*)$ analytically. Namely,
we solve the following system for $a$ and $b$:
\begin{eqnarray}\label{eq:differential-system}
\begin{cases}
\left.\frac{dR_0(y)}{dy}\right|_{y=F_0(b)}(F_0(b)-F_0(c))=R_0(F_0(b))\\
\left.\frac{dR_1(y)}{dy}\right|_{y=G_1(a)}(G_1(a)-G_1(d))=R_1(G_1(a))
\end{cases}
\end{eqnarray}
where $\left.\frac{dR_0(y)}{dy}\right|_{y=F_0(b^*)}=\beta_0^*$ and
$\left.\frac{dR_1(y)}{dy}\right|_{y=G_1(a^*)}=\beta_1^*$.
\end{remark}
Once we find $W_1(\cdot)$ and $W_0(\cdot)$, then we convert to the
original space and add back $g_1(x)$ and $g_0(x)$ respectively so
that $v_1(x)=\psi_1(x)W_1(G_1(x))+g_1(x)$ and
$v_0(x)=\varphi_0(x)W_0(F_0(x))+g_0(x)$.  Therefore, by
(\ref{eq:ost-1}) and (\ref{eq:ost-0}), the value functions
$v_1(\cdot)$ and $v_0(\cdot)$ are given by:
\begin{proposition}\label{prop:a<b}
If the optimal continuation regions for both of the value functions
are connected and if $l_c=l_d=0$, then the pair of the value
functions $v_1(x)$ and $v_0(x)$ are represented as \begin{align*}
v_1(x)=\begin{cases} \hat{v}_0(x)-H(x, 0), &x\leq a^*,\\
\hat{v}_1(x)\triangleq \psi_1(x)W_1(G_1(x))+g_1(x), &a^*<x,
\end{cases}
\end{align*}
\noindent and
\begin{align*}
v_0(x)=\begin{cases} \hat{v}_0(x)\triangleq \varphi_0(x)W_0(F_0(x))+g_0(x) &x <b^*,\\
\hat{v}_1(x)-H(x, 1), &b^*\leq x,
\end{cases}
\end{align*}
for some $a^*, b^*\in \R$ with $a^*<b^*$.
\end{proposition}
\begin{proof}
If the optimal continuation regions for both of the value functions
are connected and if $l_d=l_c=0$, then the optimal intervention
times (\ref{eq:optimal-switching-time}) have the following form:
\begin{align}\label{eq:optimal-switching-time} \theta_i^*=\begin{cases}\inf\{t>\theta_{i-1}; X_t\notin
(a^*, d)\}, &Z=1,\\
\inf\{t>\theta_{i-1}; X_t\notin (c, \hspace{0.2cm} b^*)\}, &Z=0.
\end{cases}
\end{align}
Indeed, since we have $l_c=l_d=0$, the linear majorants $W_1(\cdot)$
and $W_0(\cdot)$ pass the origins in their respective transformed
coordinates.  Hence the continuation regions shall necessarily of
the form of (\ref{eq:optimal-switching-time}).

By our construction, both $v_1(x)$ and $v_0(x)$ are continuous in
$x\in \R$.  Suppose we have $a^*> b^*$.  In this case, by the form
of the value functions, $v_0(b-)-H(b, 1, 0)=v_1(b)$.  Since the cost
function $H(\cdot)>0 $ and continuous, it follows $v_0(b-)>v_1(b)$.
On the other hand, $v_0(b+)=v_1(b)-H(b, 0, 1)$ implying
$v_0(b+)<v_1(b)$. This contradicts the continuity of $v_0(x)$. Also,
$a^*=b^*$ will lead to $v_1(x)=v_1(x)-H(x, 1, 0)$ which is
impossible. Hence if the value functions exist, then we must
necessarily have $a^*<b^*$.
\end{proof}

In relation to Proposition \ref{prop:a<b}, we have the following
observations:
\begin{remark}\normalfont
\begin{itemize}
\item[(a)]
It is obvious that
\begin{equation*}
v_0(x)=\hat{v}_0(x)>\hat{v}_0(x)-H(x, 0)=v_1(x), \quad x\in (c,
a^*),
\end{equation*}
and
\begin{equation*}
v_1(x)=\hat{v}_1(x)>\hat{v}_1(x)-H(x, 1)=v_0(x), \quad x\in (b^*,
d).
\end{equation*}
\item[(b)]
Since $u_1(x)$ is continuous in $(c, d)$, the ``obstacle"
$u_1(x)+g_1(x)-g_0(x)-H(x, 1)$ to be majorized by $u_0(x)$ on $x\in
\C_0=(c, b^*)$ is also continuous, in particular at $x=a^*$. We
proved that $u_0(x)$ always majorizes the obstacle on $(c, a^*)$.
Hence $F(a^*) \in \{y: W_0(y)>R_0(y)\}$ if there exists a linear
majorant of $R_0(y)$ in an interval of the form $(F_0(q), F_0(d))$
with some $q\in (c, d)$: otherwise, the continuity of
$u_1(x)+g_1(x)-g_0(x)-H(x, 1)$ does not hold. Similarly, we have
$F(b^*) \in \{y:W_1(y)>R_1(y)\}$ if there exists a linear majorant
of $R_0(y)$ in an interval of the form $(G_1(c), G_1(q))$.
\end{itemize}
\end{remark}

\noindent Finally, we summarize other cases than $l_c=l_d=0$:
\begin{proposition}\label{prop:odd-case}
\item[(a)] If either $l_d=+\infty$ or $l_c=+\infty$, then
$v_1(x)=v_0(x)\equiv +\infty$.
\item[(b)] If both $l_d$ and $l_c$ are finite,  then $l_d=l_c=0$.
\end{proposition}
\begin{proof}
(a) The proof is immediate by invoking Proposition \ref{prop:A5}.
(b) When $l_c$ is finite, we know by Proposition \ref{prop:A5} that
the value function $v_0(x)$ is finite.  On $x\in (c, a^*)$,
$u_1(x)+g_1(x)-g_0(x)-H(x, 1)<u_0(x)<+\infty$ is finite (see
(\ref{eq:on-(c, a)})) and thereby
\begin{equation*}
l_c=\limsup_{x\downarrow c}\frac{u_1(x)+g_1(x)-g_0(x)-H(x,
1)}{\varphi_0(x)}=0.
\end{equation*}
The same argument for $l_d=0$.
\end{proof}
\noindent Therefore, we can conclude that $l_d=0$ for the situation
where the orders of $\max(K_1(x), \psi_1(x))$ and $\psi_0(x)$ are
equal ($\Rightarrow l_d$ is finite) as described in Remark
\ref{rem:lc-ld} (a).
\section{Examples}\label{sec:example}
We recall some useful observations.  If $h(\cdot)$ is
twice-differentiable at $x\in \mathcal{I}$ and $y\triangleq F(x)$,
then we define $H(y)\triangleq h(F^{-1}(y))/\varphi(F^{-1}(y))$ and
we obtain $H^{'}(y)=m(x)$ and $H^{''}(y)=m^{'}(x)/F^{'}(x)$ with
\begin{equation}\label{eq:devH}
m(x)= \frac{1}{F^{'}(x)}\left(\frac{h}{\varphi}\right)^{'}(x), \quad
\text{and} \quad H^{''}(y) (\mathcal{A}-\alpha)h(x)\geq 0, \quad
y=F(x)
\end{equation}
with strict inequality if $H^{''}(y)\neq 0$.  These identities are
of practical use in identifying the concavities of $H(\cdot)$ when
it is hard to calculate its derivatives explicitly.  Using these
representations, we can modify (\ref{eq:differential-system}) to
\begin{eqnarray}\label{eq:modified-system}
\begin{cases}
\frac{1}{F_0'(b)}\left(\frac{r_0}{\varphi_0}\right)'(b)(F_0(b)-F_0(c))=\frac{r_0(b)}{\varphi_0(b)}\\
\frac{1}{G_1'(a)}\left(\frac{r_1}{\psi_1}\right)'(a)(G_1(a)-G_1(d))=\frac{r_1(a)}{\psi_1(a)}
\end{cases}
\end{eqnarray}
\begin{example}\normalfont \textbf{Brekke and {\O}ksendal
\cite{BO1994}}:\normalfont \hspace{0.2cm} We first illustrate our
solution method by using a resource extraction problem solved by
Brekke and {\O}ksendal~\cite{BO1994}. The price $P_t$ at time $t$
per unit of the resource follows a geometric Brownian motion. $Q_t$
denotes the stock of remaining resources in the field that decays
exponentially. Hence we have
\begin{equation*}
dP_t=\alpha P_tdt + \beta P_tdW_t\quad\text{and}\quad dQ_t=-\lambda
Q_tdt
\end{equation*}
where $\alpha, \beta,$ and $\lambda>0$ (extraction rate) are
constants. The objective of the problem is to find the optimal
switching times of resource extraction:
\begin{equation*}
v(x)=\sup_{w\in W}J^w(x)=\sup_{w\in W}\ME\left[\int_0^\infty
e^{-\rho t}(\lambda P_tQ_t-K)Z_tdt -
\sum_{i}e^{-\rho\theta_i}H(X_{\theta_i-}, Z_{\theta_i})\right]
\end{equation*}
where $rho\in\R_+$ is a discount factor with $\rho>\alpha$,
$K\in\R_+$ is the operating cost and $H(x, 0)=C\in\R_+$ and $H(x,
1)=L\in\R_+$ are constant closing and opening costs. Since $P$ and
$Q$ always show up in the form of $PQ$, we reduce the dimension by
defining $X_t=P_tQ_t$ with the dynamics:
\begin{equation*}
dX_t=(\alpha - \lambda Z_t)X_tdt + \beta X_tdW_t.
\end{equation*}
\textbf{Solution}:  
(1) We shall calculate all the necessary functions.  For $Z_t=1$
(open state), we solve $(\A_1-\rho)v(x)=0$ where
$\A_1=(\alpha-\lambda)xv'(x)+\frac{1}{2}\beta^2 x^2v''(x)$ to obtain
$\psi_1(x)=x^{\nu_+}$ and $\varphi_1(x)=x^{\nu_-}$ where $\nu_{+,
-}=\beta^{-2}\left(-\alpha+\lambda +
\frac{1}{2}\beta^2\pm\sqrt{(\alpha-\lambda-\frac{1}{2}\beta^2)^2+2\rho\beta^2}\right)$.
Similarly, for $Z_t=0$ (closed state), we solve $(\A_0-\rho)v(x)=0$
where $\A_0=\alpha xv'(x)+\frac{1}{2}\beta^2 x^2 v''(x)$ to obtain
$\psi_0(x)=x^{\mu_+}$ and $\varphi_0(x)=x^{\mu_-}$ where $\mu_{+,
-}=\beta^{-2}\left(-\alpha +
\frac{1}{2}\beta^2\pm\sqrt{(\alpha-\frac{1}{2}\beta^2)^2+2\rho\beta^2}\right)$.
Note that under the assumption $\rho>\alpha$, we have $\nu_+,
\mu_+>1$ and $\nu_-, \nu_- <0$.

By setting
$\Delta_1=\sqrt{(\alpha-\lambda-\frac{1}{2}\beta^2)^2+2\rho\beta^2}$
and $\Delta_0=\sqrt{(\alpha-\frac{1}{2}\beta^2)^2+2\rho\beta^2}$, we
have $G_1(x)=-\varphi_1(x)/\psi_1(x)=-x^{-2\Delta_1/\beta^2}$ and
$F_0(x)=\psi_0(x)/\varphi_0(x)=x^{2\Delta_0/\beta^2}$.  It follows
that $G_1^{-1}(y)=(-y)^{-\beta^2/2\Delta_1}$ and
$F_0^{-1}(y)=y^{\beta^2/2\Delta_0}$.  In this problem, we can
calculate $g_1(x), g_0(x)$ explicitly:
\begin{equation*}
g_1(x)=\ME\left[\int_0^\infty e^{-\rho s}(\lambda X_s
-K)ds\right]=\frac{x}{\rho+\lambda-\alpha}-\frac{K}{\rho}
\end{equation*}
and $g(x)=0$.  Lastly, $K_1(x)=g_0(x)-g_1(x)-H(x,
0)=-\left(\frac{x}{\rho+\lambda-\alpha}-\frac{K}{\rho}\right)-C$ and
$K_0(x)=g_1(x)-g_0(x)-H(x,
1)=\frac{x}{\rho+\lambda-\alpha}-\frac{K}{\rho}-L$.\\

\noindent (2) The state space of $X$ is $(c, d)=(0, \infty)$ and we
evaluate $l_c$ and $l_d$.  Let us first note that
$\Delta_0-\Delta_1+\lambda>0$.  Since $\lim_{x\downarrow
0}\frac{\varphi_1(x)}{\varphi_0(x)}=\lim_{x\downarrow
0}x^{\frac{\Delta_0-\Delta_1+\lambda}{\beta^2}}=0$ and
$\lim_{x\downarrow 0}(K_0(x))^+/\varphi_0(x)=0$,  we have
$l_c=l_0=0$ by (\ref{eq:sufficient-c}).  Similarly, by noting
$\lim_{x\uparrow +\infty}\frac{\psi_0(x)}{\psi_1(x)}=\lim_{x\uparrow
+\infty}x^{\frac{-(\Delta_0-\Delta_1+\lambda)}{\beta^2}}=0$ and
$\lim_{x\uparrow +\infty}(K_1(x))^+/\varphi_0(x)=0$, we have
$l_d=l_{+\infty}=0$ by (\ref{eq:sufficent-d}).\\

\noindent (3) To find the value functions together with continuation
regions, we set
\begin{align*}
r_1(x)=-\left(\frac{x}{\rho+\lambda-\alpha}-\frac{K}{\rho}\right)-C+\beta_0\psi_0(x)\quad\text{and}\quad
r_0(x)=\left(\frac{x}{\rho+\lambda-\alpha}-\frac{K}{\rho}\right)-L-\beta_1\varphi_1(x)
\end{align*}
and make transformations $R_1(y)=r_1(F^{-1}(y))/\psi_1(F^{-1}(y))$
and $R_0(y)=r_0(F^{-1}(y))/\varphi_0(F^{-1}(y))$, respectively.  We
examine the shape and behavior of the two functions $R_1(\cdot)$ and
$R_0(\cdot)$ with an aid of (\ref{eq:devH}).   By calculating $(r_0/\varphi_0)'(x)$ explicitly to examine
the derivative of $R_0(y)$,  we can find a critical point $x=q$,   at which $R_0(F(x))$ attains a local minimum and  from which  $R_0(F(x))$ is increasing monotonically on  $(F_0(q), \infty)$. Moreover, we can confirm that $\lim_{y\rightarrow
\infty}R_0'(y)=\lim_{x\rightarrow \infty}\frac{(r_0/\varphi_0)'(x)}{F_0'(x)}=0$,
which shows that there exists a finite linear majorant of $R_0(y)$.  We define
\begin{equation*}
p(x)=\beta_1 \omega
x^{\nu_-}-(\rho-\alpha)\left(\frac{x}{\rho+\lambda-\alpha}\right)+(K+\rho
L)
\end{equation*}
such that $(\A_0-\rho)r_0(x)=p(x)$ where $\omega\triangleq
\left(\rho-\frac{1}{2}\beta^2\nu_-(\nu_{-}-1)-\alpha\nu_-\right)
=\frac{1}{2\beta^2}(\Delta_0-\Delta_1+\lambda)(\Delta_0+\Delta_1-\lambda)>0$.
By the second identity in (\ref{eq:devH}), the sign of the second
derivative $R_0''(y)$ is the same as the sign of $p(x)$. It is easy
to see that $p(x)$ has only one critical point.  For any
$\beta_1<0$, the first term is dominant as $x\rightarrow 0$, so that
$\lim_{x\downarrow 0} p(x)<0$.  As $x$ gets larger, for $|\beta_1|$
sufficiently small, $p(x)$ can take positive values, providing two
positive roots, say $x=k_1, k_2$ with $k_1<k_2$. We also have
$\lim_{x\rightarrow +\infty}p(x)=-\infty$. In this case, $R_0(y)$ is
concave on $(0, F(k_1)\cup (F(k_2), +\infty)$ and convex on
$(F(k_1), F(k_2))$.  Since we know that $R_0(y)$ attains a local
minimum at $y=F(q)$, we have $q<k_2$, and it implies that there is
one and only on tangency point of the linear majorant $W(y)$ and
$R_0(y)$ on $(F(q), \infty)$, so that the continuation region is of
the form $(0, b^*)$.

>From this analysis of the derivatives of $R_0(y)$, there is only one
tangency point of the linear majorant $W_0(y)$ and $R_0(y)$. (See
Figure \ref{fig:1}-(a)). A similar analysis shows that there is only
one tangency point of the linear majorant $W_1(y)$ and $R_1(y)$.
(See Figure \ref{fig:1}-(b)).
\begin{figure}[h]\label{fig:1}
\begin{center}
\begin{minipage}{0.45\textwidth}
\centering \includegraphics[scale=0.75]{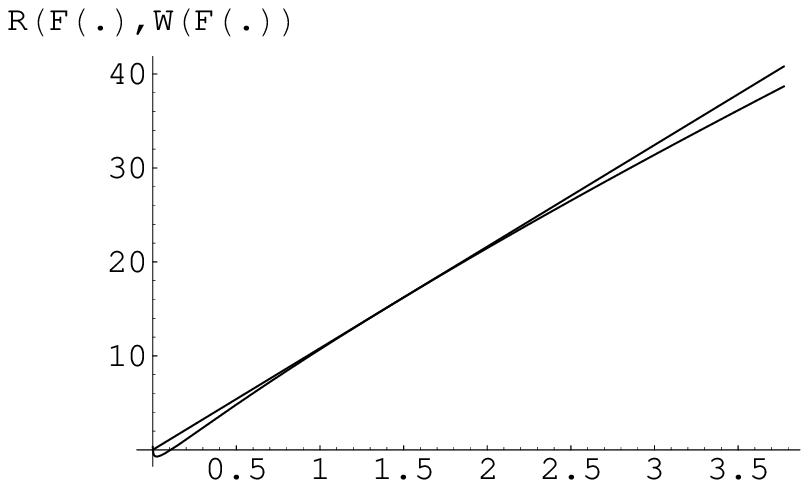} \\
(a)
\end{minipage}
\begin{minipage}{0.45\textwidth}
\centering \includegraphics[scale=0.75]{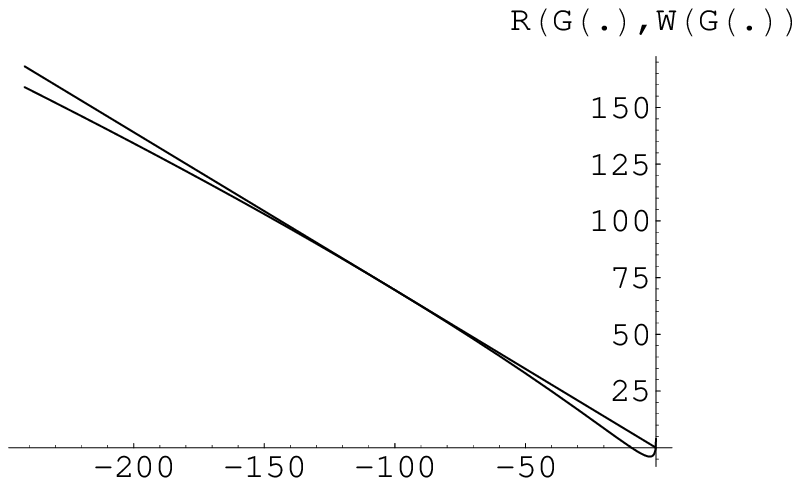} \\
(b)
\end{minipage}
\begin{minipage}{0.45\textwidth}
\centering \includegraphics[scale=0.75]{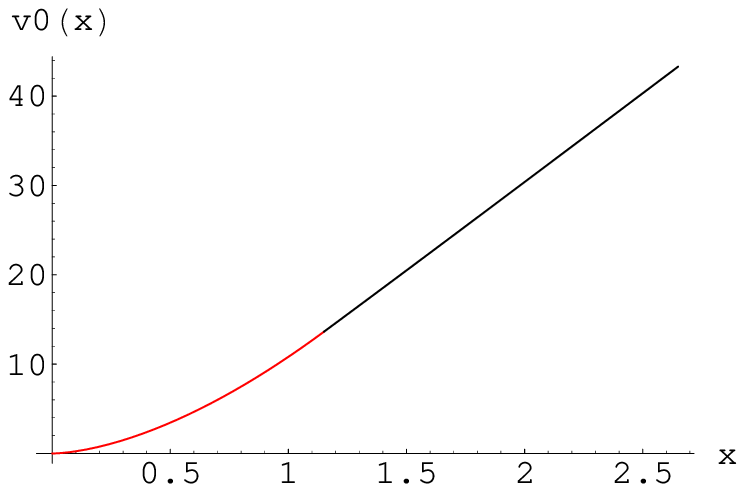} \\
(c)
\end{minipage}
\begin{minipage}{0.45\textwidth}
\centering \includegraphics[scale=0.75]{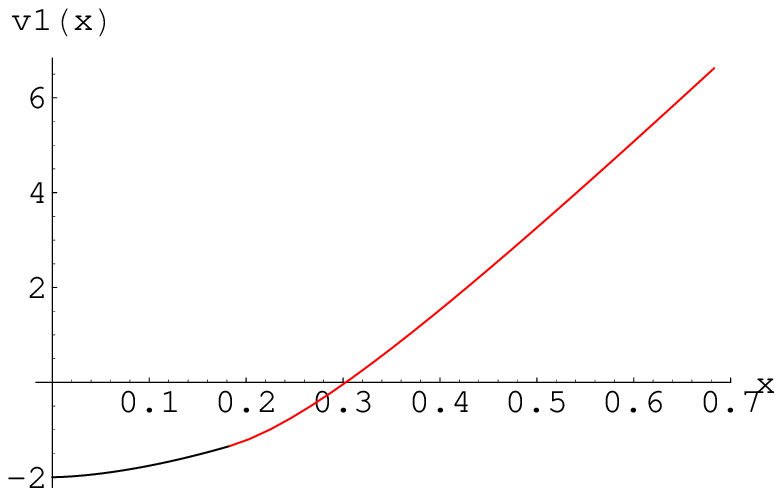} \\
(d)
\end{minipage}
\caption{\small A numerical example of resource extraction problem.
with parameters $(\alpha, \beta, \lambda, \rho, K, L, C)=(0.01,
0.25, 0.01, 0.05, 0.4, 2, 2)$(a) The smallest linear majorant
$W_0(F_0(x))$ and $R_0(F_0(x))$ with $b^*=1.15042$ and
$\beta_0^*=10.8125$.  (b)The smallest linear majorant $W_1(G_1(x))$
and $R_1(G_1(x))$ with $a^*=0.18300$ and $\beta_1^*=-0.695324$.  (c)
The value function $v_0(x)$. (d) The value function $v_1(x)$.}
\end{center}
\end{figure}\\
\noindent (4) By solving the system of equations
(\ref{eq:differential-system}), we can find $(a^*, b^*, \beta_0^*,
\beta_1^*)$.  We transform back to the original space to find
\begin{align*}
\hat{v}_1(x)&=\psi_1(x)W_1(G_1(x))+g_1(x)=\psi_1(x)\beta_1^*G_1(x)+g_1(x)\\
&=-\beta_1^*\varphi_1(x)+g_1(x)=-\beta_1^*
x^{\nu_-}+\left(\frac{x}{\rho+\lambda-\alpha}-\frac{K}{\rho}\right),
\end{align*}
and
\begin{align*}
\hat{v}_0(x)&=\varphi_0(x)W_0(F_0(x))+g_0(x)=\varphi_0(x)\beta_0^*F_0(x)+g_0(x)=\beta_0^*\psi_0(x)+g_0(x)=\beta_0^*x^{\mu_{+}}.
\end{align*}
Hence the solution is
\begin{align} \nonumber
\begin{aligned}
    v_1(x) &=  \begin{cases}
                 \beta_0^*x^{\mu_+}-C, & x\leq a^*, \\
                 -\beta_1^*
x^{\nu_-}+\left(\frac{x}{\rho+\lambda-\alpha}-\frac{K}{\rho}\right),
                 &x>a^*,\end{cases}
\hspace{0.4cm}
    v_0(x) &= \begin{cases}
                 \beta_0^*x^{\mu_{+}}, & x\leq b^*, \\
                 -\beta_1^*
x^{\nu_-}+\left(\frac{x}{\rho+\lambda-\alpha}-\frac{K}{\rho}\right)-L,
                 &x>b^*,
     \end{cases}
\end{aligned}
\end{align}
which agrees with Brekke and {\O}kesendal \cite{BO1994}.
\end{example}
\begin{example}\normalfont \textbf{Ornstein-Uhrenbeck process}:
\hspace{0.2cm}  We shall consider a new problem involving an
Ornstein-Uhrenbeck process. Consider a firm whose revenue solely
depends on the price of one product. Due to its cyclical nature of
the prices, the firm does not want to have a large production
facilty and decides to rent additional production facility when the
price is favorable. The revenue process to the firm is
\begin{equation*} dX_t=\delta(m-X_t-\lambda Z_t)dt+\sigma dW_t,
\end{equation*}
where $\lambda=r/\delta$ with $r$ being a rent per unit of time. The
firm's objective is to maximize the \emph{incremental} revenue
generated by renting the facility until the time $\tau_0$ when the
price is at an intolerably low level. Without loss of generality, we
set $\tau_0=\inf\{t>0: X_t=0\}$.  We keep assuming constant
operating cost $K$, opening cost, $L$ and closing cost $C$. Now the
value function is defined as
\begin{equation*}
v(x)=\sup_{w\in W}J^w(x)=\sup_{w\in W}\ME\left[\int_0^{\tau_0}
e^{-\alpha t}(X_t-K)Z_tdt -
\sum_{\theta_i<\tau_0}e^{-\alpha\theta_i}H(X_{\theta_i-},
Z_{\theta_i})\right].
\end{equation*}
\textbf{Solution}:  
(1) We denote, by $\tilde{\psi}(\cdot)$ and
$\tilde{\varphi}(\cdot)$, the functions of the fundamental solutions
for the auxiliary process $P_t\triangleq(X_t-m+\lambda)/\sigma,
t\geq 0$, which satisfies $dP_t=-\delta P_t dt+dW_t$.  For every
$x\in \mathbb{R}$,
\begin{equation*}
\tilde{\psi}(x)=e^{\delta
x^2/2}\mathcal{D}_{-\alpha/\delta}(-x\sqrt{2\delta})\quad\text{and}\quad
\tilde{\varphi}(x)=e^{\delta
x^2/2}\mathcal{D}_{-\alpha/\delta}(x\sqrt{2\delta}),
\end{equation*}
which leads to $\psi_1(x)=\tilde{\psi}((x-m+\lambda)/\sigma)$,
$\varphi_1(x)=\tilde{\varphi}((x-m+\lambda)/\sigma)$,
$\psi_0(x)=\tilde{\psi}((x-m)/\sigma)$, and
$\varphi_0(x)=\tilde{\varphi}((x-m)/\sigma)$ where
$\mathcal{D}_\nu(\cdot)$ is the parabolic cylinder function; (see
Borodin and Salminen (2002, Appendices 1.24 and 2.9) and Carmona and
Dayanik (2003, Section 6.3)).  By using the relation
\begin{equation} \label{eq:Dft}
\mathcal{D}_\nu(z)=2^{-\nu/2}e^{-z^2/4}\mathcal{H}_\nu(z/\sqrt{2}),
\quad z\in\mathbb{R}
\end{equation}
in terms of the Hermite function $\mathcal{H}_\nu$ of degree $\nu$
and its integral representation
\begin{equation}\label{eq:Hermite}
\mathcal{H}_\nu(z)=\frac{1}{\Gamma(-\nu)}\int_0^\infty
e^{-t^2-2tz}t^{-\nu-1}dt, \quad \text{Re}(\nu)<0,
\end{equation}
(see for example, Lebedev(1972, pp 284, 290)). Since
$\ME[X_t]=e^{-\delta t}x+(1-e^{-\delta t})(m-\lambda)$, we have
$g_0(x)=0$ and
$g_1(x)=\frac{x-(m-\lambda)}{\delta+\alpha}+\frac{m-\lambda-K}{\alpha}$.\\

\noindent(2) The state space of $X$ is $(c, d)=(0, +\infty)$.  Since
the left boundary $0$ is the absorbing, the linear majorant passes
$(0, F_0(0))$.  Since $\lim_{x\rightarrow
+\infty}\psi_0(x)/\psi_1(x)=0$, we have
$l_d=0$.\\

\noindent (3) We formulate
\begin{equation*}
r_1(x)=-\left(\frac{x-(m-\lambda)}{\delta+\alpha}+\frac{m-\lambda-K}{\alpha}\right)-C+\beta_0\psi_0(x)
\end{equation*}
and 
\begin{equation*}
r_0(x)=\left(\frac{x-(m-\lambda)}{\delta+\alpha}+\frac{m-\lambda-K}{\alpha}\right)-L-\beta_1\varphi_1(x)
\end{equation*}
and make transformations: $R_1(y)=r_1(F^{-1}(y))/\psi_1(F^{-1}(y))$
and $R_0(y)=r_0(F^{-1}(y))/\varphi_0(F^{-1}(y))$, respectively.  We
examine the shape and behavior of the two functions $R_1(\cdot)$ and
$R_0(\cdot)$ with an aid of (\ref{eq:devH}).  First we check the
sign of $R_0'(y)$ and  find a critical point $x=q$,   at which
$R_0(F(x))$ attains a local minimum and  from which  $R_0(F(x))$ is
increasing monotonically on  $(F_0(q), \infty)$.  It can be shown
that $R_0^{'}(+\infty)=0$ by using (\ref{eq:Dft}) and
(\ref{eq:Hermite}) and the identity
$\mathcal{H}'_\nu(z)=2\nu\mathcal{H}_{\nu-1}(z), z\in\mathbb{R}$
(see Lebedev (1972, p.289), for example.) This shows that there must
exist a (finite) linear majorant of $R_0(y)$ on $(F(q), \infty)$.
To check convexity of $R_0(y)$, we define
\begin{equation*}
p(x)=-\frac{\sigma^2\beta_1}{2}\varphi_1''(x)+\delta(m-x-\lambda)\left(\frac{1}{\delta+\alpha}-\beta_1\varphi_1'(x)\right)
-\alpha r_0(x)
\end{equation*}
such that $(\A_0-\alpha)r_0(x)=p(x)$.  We can show easily
$\lim_{x\rightarrow +\infty}p(x)=-\infty$ since
$\varphi_1(+\infty)=\varphi_1'(+\infty)=\varphi_1''(+\infty)=0$. Due
to the monotonicity of $\varphi_1(x)$ and its derivatives, $p(x)$
can have at most one critical point and $p(x)=0$ can have one or two
positive roots depending on the value of $\beta_1$. In either case,
let us call the largest positive root $x=k_2$. We also have
$\lim_{x\rightarrow +\infty}p(x)=-\infty$. 
Since we know that $R_0(y)$ attains a local minimum at $y=F(q)$ and
is increasing thereafter, we have $q<k_2$.  It follows that there is
one and only on tangency point of the linear majorant $W(y)$ and
$R_0(y)$ on $(F(q), \infty)$, so that the continuation region is of
the form $(0, b^*)$. A similar analysis shows that there is only
one tangency point of the linear majorant $W_1(y)$ and $R_1(y)$.\\

\noindent (4)  Solving (\ref{eq:modified-system}), we we can find
$(a^*, b^*, \beta_0^*, \beta_1^*)$.  We transform back to the
original space to find
\begin{align*}
\hat{v}_1(x)&=\psi_1(x)W_1(G_1(x))+g_1(x)=\psi_1(x)\beta_1^*G_1(x)+g_1(x)=-\beta_1^*\varphi_1(x)+g_1(x)\\
&=-\beta_1^*
e^{\frac{\delta(x-m+\lambda)^2}{2\sigma^2}}\mathcal{D}_{-\alpha/\delta}\left(\frac{(x-m+\lambda)\sqrt{2\delta}}{\sigma}\right)
+\frac{x-(m-\lambda)}{\delta+\alpha}+\frac{m-\lambda}{\alpha}
\end{align*}
and
\begin{align*}
\hat{v}_0(x)&=\varphi_0(x)W_0(F_0(x))+g_0(x)=\varphi_0(x)\beta_0^*(F_0(x)-F_0(0))+g_0(x)\\
&=\beta_0^*\{\psi_0(x)-F_0(0)\varphi_0(x)\}+g_0(x)\\
&=\beta_0^*e^{\frac{\delta}{2}\frac{(x-m+\lambda)^2}{\sigma^2}}\left\{\mathcal{D}_{-\alpha/\delta}\left(-\left(\frac{x-m+\lambda}{\sigma}\right)\sqrt{2\delta}\right)
-F(0)\mathcal{D}_{-\alpha/\delta}\left(\left(\frac{x-m}{\sigma}\right)\sqrt{2\delta}\right)\right\}.
\end{align*}
Hence the solution is, using the above functions,
\begin{align} \nonumber
\begin{aligned}
    v_1(x) &=  \begin{cases}
                 \hat{v}_0(x)-C, & x\leq a^*, \\
                 \hat{v}_1(x),
                 &x>a^*,\end{cases}
\hspace{0.4cm}
    v_0(x) &= \begin{cases}
                 \hat{v}_0(x), & x\leq b^*, \\
                 \hat{v}_1(x)-L,
                 &x>b^*.
     \end{cases}
\end{aligned}
\end{align}
See Figure \ref{fig:2} for a numerical example.
\begin{figure}[h]\label{fig:2}
\begin{center}
\begin{minipage}{0.45\textwidth}
\centering \includegraphics[scale=0.75]{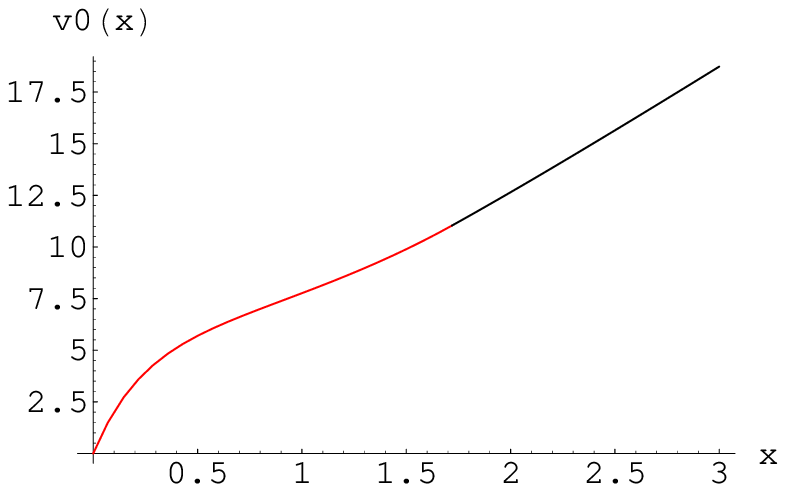} \\
(a)
\end{minipage}
\begin{minipage}{0.45\textwidth}
\centering \includegraphics[scale=0.75]{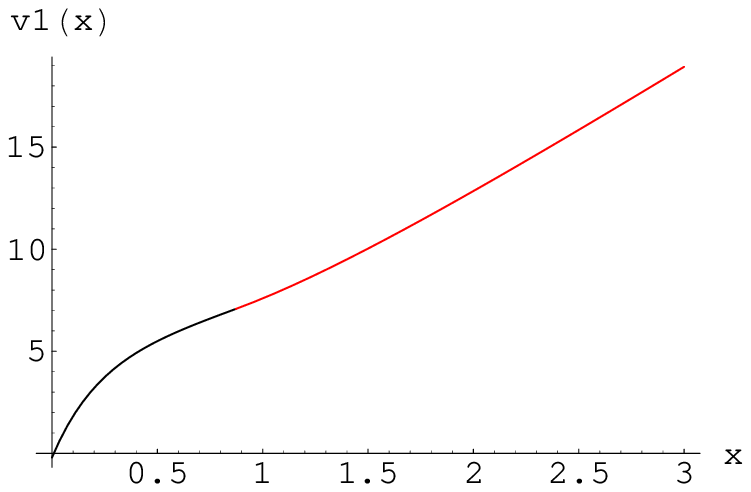} \\
(b)
\end{minipage}
\caption{\small A numerical example of leasing production facility
problem with parameters $(m, \alpha, \sigma, \delta, \lambda,  K, L,
C)=(5, 0.105, 0.35, 0.05, 4, 0.4, 0.2, 0.2)$: (a) The value function
$v_0(x)$ with $b^*=1.66182$ and $\beta_0^*=144.313$.  (b)The value
function $v_1(x)$ with $a^*=0.781797$ and $\beta_1^*=-2.16941$.}
\end{center}
\end{figure}
\end{example}
\section{Extensions and conclusions}\label{sec:last-section}
\subsection{An extension to the case of $k\geq 2$}
 It is not difficult to extend to a general case of
$k\geq 2$ where more than one switching opportunities are available.
But we put a condition that $z\in\mathcal{Z}$ is of the form
$z=(a_1, a_2,...., a_k)$ where only one element of this vector is
$1$ with the rest being zero, i.e., $z=(0, 0, 0,...., 1, 0, 0)$ for
example.

We should introduce the switching operator $\mathcal{M}_0$ on $h\in
\mathcal{H}$,
\begin{equation}\label{eq:switch-operator}
\mathcal{M}_0h(u,z)=\max_{\zeta\in \mathcal{Z}\setminus
\{z\}}\left\{h(u, \zeta)-H(u, z; \zeta)\right\}.
\end{equation}
In words, this operator would calculate which production mode should
be chosen by moving from the current production mode $z$.   Now the
recursive optimal stopping (\ref{eq:w-function}) becomes
\begin{equation*}
w_{n+1}(x)\triangleq\sup_{\tau\in \mathcal{S}}\ME\left[\int_0^\tau
e^{-\alpha s}f(X_s)ds + e^{-\alpha\tau}
\mathcal{M}w_{n}(X_{\tau})\right].
\end{equation*}
Accordingly, the optimization procedure will become two-stage. To
illustrate this, we suppose $k=2$ so that $i=0, 1,$ and $2$. By
eliminating the integral in (\ref{eq:switch-operator}), we redefine
the switching operator,
\begin{equation}
\M h_z(x)\triangleq \max_{\zeta\in\mathcal{Z}\setminus
\{z\}}\left\{h_\zeta(x)+g_{\zeta}(x)-g_z(x)-H(x, z, \zeta)\right\},
\end{equation}
where
\begin{equation*}
g_z(x)\triangleq \sup_{w\in W_0}J_z^w(x)=\ME\left[\int_0^\infty
e^{-\alpha s}f(X_s)ds \right].
\end{equation*}
Hence (\ref{eq:u-1}) will be modified to
$u_z(x)=\ME[e^{-\alpha\tau}\M u_z(X_\tau)].$
It follows that our system of equations (\ref{eq:system}) is now
\begin{align}
\begin{cases}\label{eq:3-system}
\bar{v}_2(x)\triangleq \sup_{\tau\in\S}\ME\left[e^{-\alpha\tau}
\M\bar{v_2}(X_{\tau})\right]\\
\bar{v}_1(x)\triangleq \sup_{\tau\in\S}\ME\left[e^{-\alpha\tau}
\M\bar{v_1}(X_{\tau})\right]\\
\bar{v}_0(x)\triangleq \sup_{\tau\in\S}\ME\left[e^{-\alpha\tau}
\M\bar{v_0}(X_{\tau})\right]\\
\end{cases}
\end{align}
The first stage is optimal stopping problem.  One possibility of
switching production modes is $(0\rightarrow 1, 1\rightarrow 2,
2\rightarrow 0)$. First, we fix this switching scheme, say $c$, and
solve the system of equations (\ref{eq:3-system}) as three optimal
stopping problems. All the arguments in Section 2.3 hold. This
first-stage optimization will give $(x_0^*(c), x_1^*(c), x_2^*(c),
\beta_0^*(c), \beta_1^*(c), \beta_2^*(c))$, where $x_i$'s are
switching boundaries, depending on this switching scheme $c$.

Now we move to another switching scheme $c'$ and solve the system of
optimal stopping problems until we find the optimal scheme.

\subsection{Conclusions}
We have studied optimal switching problems for one-dimensional
diffusions.  We characterize the value function as linear functions
in their respective spaces, and provide a direct method to find the
value functions and the opening and switching boundaries at the same
time.  Using the techniques we developed here as well as the ones in
Dayanik and Karazas \cite{DK2003} and Dayanik and Egami
\cite{DE2005}, we solved two specific problems, one of which
involves a mean-reverting process.  This problem might be hard to
solve with just the HJB equation and the related quasi-variational
inequalities.  Finally, an extension to more general cases is
suggested.  We believe that this direct method and the new
characterization will expand the coverage of solvable problems in
the financial engineering and economic analysis.

\begin{appendix}
\section{Summary of Optimal Stopping Theory}\label{appx:ost}
Let $(\Omega, \F, \p)$ be a complete probability space with a
standard Brownian motion $W=\{W_t; t\geq 0\}$ and consider the
diffusion process $X^0$ with state pace $\mathcal{I}\subseteq
\mathbb{R}$ and dynamics
\begin{equation}\label{eq:process}
dX^0_t=\mu(X^0_t)dt + \sigma(X^0_t)dW_t
\end{equation}
for some Borel functions $\mu :\mathcal{I}\rightarrow \mathbb{R}$
and $\sigma :\mathcal{I}\rightarrow (0, \infty)$.  We emphasize here
that $X^0$ is an uncontrolled process.  We assume that $\mathcal{I}$
is an interval with endpoints $-\infty\leq a < b \leq+\infty$, and
that $X^0$ is regular in $(a, b)$; in other words, $X^0$ reaches $y$
with positive probability starting at $x$ for every $x$ and $y$ in
$(a,b)$.  We shall denote by $\mathbb{F}=\{\mathcal{F}_t\}$ the
natural filtration generated by $X^0$.

Let $\alpha \geq 0$ be a real constant and $h(\cdot)$ a Borel
function such that $\ME[e^{-\alpha \tau}h(X^0_{\tau})]$ is
well-defined for every $\mathbb{F}$-stopping time $\tau$ and $x\in
\mathcal{I}$.  Let $\tau_y$ be the first hitting time of $y\in
\mathcal{I}$ by $X^0$, and let $c\in\mathcal{I}$ be a fixed point of
the state space.  We set:
\begin{align} \nonumber
\begin{aligned}
    \psi(x) &=  \begin{cases}
                 \ME[e^{-\alpha\tau_c}1_{\{\tau_c<\infty\}}], & x\leq c, \\
                 1/\E^{c}[e^{-\alpha\tau_x}1_{\{\tau_x<\infty\}}],
                 &x>c,\end{cases}
\hspace{0.4cm}
    \varphi(x) &= \begin{cases}
                 1/\E^{c}\left[e^{-\alpha\tau_x}1_{\{\tau_x<\infty\}}\right], & x\leq c, \\
                 \ME[e^{-\alpha\tau_c}1_{\{\tau_c<\infty\}}],
                 &x>c,
     \end{cases}
\end{aligned}
\end{align}
and
\begin{align} \label{eq:F}
F(x)&\triangleq\frac{\psi(x)}{\varphi(x)}, \hspace{0.5cm} x\in
\mathcal{I}.
\end{align}
Then $F(\cdot)$ is continuous and strictly increasing.  It should be
noted that $\psi(\cdot)$ and $\varphi(\cdot)$ consist of an
increasing and a decreasing solution of the second-order
differential equation $(\mathcal{A}-\alpha)u=0$ in $\mathcal{I}$
where $\mathcal{A}$ is the infinitesimal generator of $X^0$.  They
are linearly independent positive solutions and uniquely determined
up to multiplication.  For the complete characterization of
$\psi(\cdot)$ and $\varphi(\cdot)$ corresponding to various types of
boundary behavior, refer to It\^{o} and McKean \cite{IM1974}.

Let $F :[c, d]\rightarrow\mathbb{R}$ be a strictly increasing
function.  A real valued function $u$ is called \emph{$F$-concave}
on $[c, d]$ if, for every $a\leq l<r\leq b$ and $x\in[l, r]$,
\begin{equation*}
u(x)\geq
u(l)\frac{F(r)-F(x)}{F(r)-F(l)}+u(r)\frac{F(x)-F(l)}{F(r)-F(l)}.
\end{equation*}
We denote by
\begin{equation}\label{eq:value}
V(x)\triangleq \sup_{\tau\in\S}\ME[\D h(X^0_\tau)], \hspace{0.5cm}
x\in[c, d]
\end{equation}
the value function of the optimal stopping problem with the reward
function $h(\cdot)$ where the supremum is taken over the class $\S$
of all $\mathbb{F}$-stopping times.  Then we have the following
results, the proofs of which we refer to Dayanik and
Karatzas~\cite{DK2003}.
\begin{proposition}\normalfont\label{prop:A1}
For a given function $U$: $[c, d]\rightarrow[0,+\infty)$ the
quotient $U(\cdot)/\varphi(\cdot)$ is an $F$-concave function if and
only if $U(\cdot)$ is $\alpha$-excessive, i.e.,
\begin{align}
U(x)\geq \ME[e^{-\alpha\tau}U(X^0_\tau)], &\forall \tau\in\S,
\forall x\in[c,d].
\end{align}
\end{proposition}
\begin{proposition}\normalfont \label{prop:A2}
The value function $V(\cdot)$ of (\ref{eq:value}) is the smallest
nonnegative majorant of $h(\cdot)$ such that
$V(\cdot)/\varphi(\cdot)$ is $F$-concave on $[c,d]$.
\end{proposition}
\begin{proposition}\normalfont \label{prop:A3}
Let $W(\cdot)$ be the smallest nonnegative concave majorant of
$H\triangleq (h/\varphi)\circ F^{-1}$ on $[F(c), F(d)]$, where
$F^{-1}(\cdot)$ is the inverse of the strictly increasing function
$F(\cdot)$ in (\ref{eq:F}).  Then $V(x)=\varphi(x)W(F(x))$ for every
$x\in[c, d]$.
\end{proposition}
\begin{proposition}\normalfont \label{prop:A4}
Define
\begin{equation}\label{eq:opt}
S\triangleq\{x\in[c, d]: V(x)=h(x)\}, \hspace{0.5cm}\text{and}
\hspace{0.5cm} \tau^{*}\triangleq \inf\{t\geqq0: X^0_t\in S\}.
\end{equation}
If $h(\cdot)$ is continuous on $[c, d]$, then $\tau^{*}$ is an
optimal stopping rule.
\end{proposition}
When both boundaries are natural, we have the following results:
\begin{proposition} \normalfont \label{prop:A5}
We have either $V\equiv 0$ in $(c, d)$ or $V(x)<+\infty$ for all
$(c, d)$.  Moreover, $V(x)<+\infty$ for every $x\in (c, d)$ if and
only if
\begin{equation}\label{eq:test}
l_c\triangleq \limsup_{x\downarrow c}\frac{h^+(x)}{\varphi(x)}\quad
\text{and}\quad l_d\triangleq \limsup_{x\uparrow
d}\frac{h^+(x)}{\psi(x)}
\end{equation}
are both finite.
\end{proposition}
In the finite case, furthermore,
\begin{proposition}\normalfont \label{prop:A6}
The value function $V(\cdot)$ is continuous on $(c, d)$.  If $h: (c,
d)\rightarrow \R$ is continuous and $l_c=l_d=0$, then $\tau^*$ of
(\ref{eq:opt}) is an optimal stopping time.
\end{proposition}
\begin{proposition}\normalfont\label{prop:A7}
Suppose that $l_c$ and $l_d$ are finite and one of them is strictly
positive, and $h(\cdot)$ is continuous.  Define the continuation
region $C\triangleq (c, d)\setminus\Gamma$.  Then $\tau^*$ of
(\ref{eq:opt}) is an optimal stopping time, if and only if
\begin{align*}
&\text{there is no $r\in (c, d)$ such that $(c, r)\subset C$ if
$l_c>0$}\quad \text{and}\\
&\text{there is no $l\in (c, d)$ such that $(l, d)\subset C$ if
$l_d>0$}.
\end{align*}
\end{proposition}

\end{appendix}
\bibliography{switch}
\end{document}